\font \sevenrm=cmr7
\font \eightrm=cmr8

\font \eightbf=cmbx8
\font \bigbf=cmbx10 scaled \magstep1
\font \Bigbf=cmbx10 scaled \magstep2

%
\font \tengoth=eufm10
\font \sevengoth=eufm7
\font \fivegoth=eufm5

\newfam\gothfam
\textfont \gothfam=\tengoth
\scriptfont \gothfam=\sevengoth
\scriptscriptfont \gothfam=\fivegoth

%
\font \tenmath=msbm10
\font \sevenmath=msbm7
\font \fivemath=msbm5

\newfam\mathfam
\textfont \mathfam=\tenmath
\scriptfont \mathfam=\sevenmath
\scriptscriptfont \mathfam=\fivemath
\def\math{\fam\mathfam\tenmath}
%
%
%
%
\def\titre#1{\centerline{\Bigbf #1}\nobreak\nobreak\vglue 10mm\nobreak}

\def\paragraphe#1{\bigskip\goodbreak {\bigbf #1}\nobreak\vglue 12pt\nobreak}

\def\ssq{\smallskip\qquad}

%
%
\def\th#1{\bigskip\goodbreak {\bf Theorem #1.} \par\nobreak \sl }
\def\prop#1{\bigskip\goodbreak {\bf Proposition #1.} \par\nobreak \sl }
\def\lemme#1{\bigskip\goodbreak {\bf Lemma #1.} \par\nobreak \sl }
\def\cor#1{\bigskip\goodbreak {\bf Corollary #1.} \par\nobreak \sl }
\def\dem{\bigskip\goodbreak \it Proof. \rm}
\def\ndem{\bigskip\goodbreak \rm}
\def\qed{\par\nobreak\hfill $\bullet$ \par\goodbreak}
%
%
\def\uple#1#2{#1_1,\ldots ,{#1}_{#2}}
\def \restr#1{\mathstrut_{\textstyle |}\raise-6pt\hbox{$\scriptstyle #1$}}
\def \srestr#1{\mathstrut_{\scriptstyle |}\hbox to -1.5pt{}\raise-4pt\hbox{$\scriptscriptstyle #1$}}
\def \inver{^{-1}}
\def\surj#1{\mathop{\hbox to #1 mm{\rightarrowfill\hskip 2pt\llap{$\rightarrow$}}}\limits}
\def\frac#1#2{{\textstyle {#1\over #2}}}
\def\R{{\math R}}

\def\Z{{\math Z}}

\def\fleche#1{\mathop{\hbox to #1 mm{\rightarrowfill}}\limits}
%
%
\def \g#1{\hbox{\tengoth #1}}
\def \sg#1{\hbox{\sevengoth #1}}

\def\Cal #1{{\cal #1}}
%
%

\def \smop#1{\mathop{\hbox{\sevenrm #1}}\nolimits}

%
%
\def \bib #1{\null\medskip \strut\llap{#1.\quad}}
\def\cite#1{[#1]}
\magnification=\magstep1
\parindent=0cm

\input xy
\xyoption{all}
\def\dm{{1\over 2}}
\def\hdm{{\hbar\over 2}}
\def\diagramme #1{\vskip 4mm \centerline {#1} \vskip 4mm}

\titre{Poisson bracket, deformed bracket}
\vskip -5mm
\titre{and gauge group actions}
\vskip -5mm
\titre{in Kontsevich deformation quantization}
\bigskip
\centerline{Dominique Manchon \footnote * {\eightrm CNRS - Institut Elie Cartan, BP 239, F54506 Vandoeuvre CEDEX. manchon@iecn.u-nancy.fr}}
\vskip 25mm
{\baselineskip=10pt\eightbf Abstract \eightrm : We express the difference between Poisson bracket and deformed bracket for Kontsevich deformation quantization on any Poisson manifold by means of second derivative of the formality quasi-isomorphism. The counterpart on star products of the action of formal diffeomorphisms on Poisson formal bivector fields is also investigated.
\smallskip
\eightbf Key words : \eightrm Poisson manifold, deformation quantization, star product, formality, gauge transformation, super-grouplike element.
\par}
\vskip 12mm
\paragraphe{Introduction}
\qquad The existence of a star product on any Poisson manifold $(M,\gamma)$ is derived from the more general formality theorem of M. Kontsevich \cite K, which stipulates the existence of a $L_\infty$-quasi-isomorphism $\Cal U$ (cf. \S\ II) from the differential graded Lie algebra of polyvector fields on any manifold $M$ (with vanishing differential an Schouten bracket) into the differential graded Lie algebra of polydifferential operators on $M$ (with Hochschild differential and Gerstenhaber bracket).
\ssq
Given such a $L_\infty$-quasi-isomorphism there is a canonical and explicit way to produce a star product $*=*_\gamma$ from Poisson bivector field, and more generally from any formal Poisson bivector field $\gamma$. We briefly recall this construction on \S ~II. We call {\sl formality star product\/} any star product on $M$ obtained that way. Due to the fact that $\Cal U$ is a $L_\infty$-quasi-isomorphism, any star product is gauge-equivalent to a formality star product \cite {K \S ~4.4}, \cite{AMM \S~A.2}.
\ssq
Let $f,g\in C^\infty(M)[[\hbar]]$, and let $H_f=[\gamma,f],\ H_g=[\gamma,g]$ the associated hamiltonian (formal) vector fields. We compute here the second derivative at $\hbar \gamma$ of quasi-isomorphism $\Cal U$ evaluated at $(H_f,H_g)$, and more generally at $(Y,H_g)$ where $Y$ stands for any formal vector field on $M$. The main result is theorem III.3, consisting of three equations, which in turn imply the following formula, relating Poisson bracket, deformed bracket, tangent map $\Phi$ at $\hbar \gamma$ of $\Cal U$ and second derivative $\Psi$ at $\hbar \gamma$ of $\Cal U$~:
$$\Psi(H_f.H_g)={1\over \hbar}\Bigl(\Phi(\{f,g\})-{\Phi(f)*\Phi(g)-\Phi(g)*\Phi(f)\over\hbar}\Bigr).$$
\qquad
There is another consequence of theorem III.3 in terms of gauge group action~: namely we try to understand the star product $*_{g.\gamma}$ obtained from formal Poisson bivector field $g.\gamma$ where $g$ is a formal diffeomorphism of the manifold $M$. Formal diffeomorphisms also act naturally on star-products via action on $C^\infty(M)$, but it is quite obvious that $*_{g.\gamma}$ is not the image of $*_\gamma$ by the action of $g$ in that sense.
\ssq
Gauge group $G_1$ of formal diffeomorphisms however acts on formality star products in a more subtle way~: considering the set of all formality star-products as a formal pointed manifold $(FSP)$, the action we seek amounts to a nontrivial embedding of $G_1$ into the group $\Cal G$ of formal diffeomorphisms of $(FSP)$.

\paragraphe{I. Super-grouplike elements in cofree cocommutative graded coalgebras}
\qquad Let $V=\bigoplus V^{(n)}$ be a $\Z$-graded vector space over a field $k$ with zero characteristic, and let $\Cal C=S(V)$ its symmetric algebra in the graded sense. We will throughout the paper denote by $\pi$ the projection of $S(V)$ onto $V$. Defining a coproduct on elements of $V$ by~:
$$\Delta(x)=x\otimes 1+1\otimes x$$
and extending it to an algebra morphism from $S(V)$ to the tensor product $S(V)\otimes S(V)$ (in the graded sense) we endow $S(V)$ with a structure of graded bialgebra. The set of primitive elements is precisely $V$, and the co-unity is given by the projection on constants.
\ssq
Let $\g m$ be a (projective limit of) commutative finite dimensional nilpotent algebra(s). We will consider the (completed) tensor product $V\widehat\otimes\g m$ as a topologically free $\g m$-module and we will see the topologically free $\g m$-module $\Cal C_{\sg m}=S(V)\widehat\otimes \g m\oplus k.1$ as a topological bialgebra over $\g m$. 
\ssq Let $\tau:\Cal C\otimes \Cal C\rightarrow \Cal C\otimes \Cal C$ be the signed flip defined by~:
$$\tau(v\otimes w)=(-1)^{|v||w|}w\otimes v.$$
A nonzero element $v\in\Cal C_{\g m}$ will be called {\sl super-grouplike} if we have~:
$$\Delta v={I+\tau\over 2}(v\otimes v).$$
As an example, any even grouplike element is super-grouplike, as well as any $1+x$ where $x\in V\widehat\otimes \g m$ and $x$ odd.
\prop{I.1}
Any super-grouplike element in the coalgebra $\Cal C_{\sg m}$ is of the form~:
$$g=e^{.v}=1+v+\cdots +{1\over n!}(v\cdots v)_{\hbox{\sevenrm n times}}+\cdots$$
with $v\in V\widehat\otimes \g m$,
and conversely any such exponential is super-grouplike.
\dem
Consider the decomposition $g=g_++g_-$ of our supergrouplike element into its even and odd components. We have then~:
$${1+\tau\over 2}\Delta(g)={1+\tau\over 2}\Delta(g_++g_-)=g_+\otimes g_++g_-\otimes g_++g_+\otimes g_-,$$
hence~:
$$\eqalign{\Delta (g_+)         &=g_+\otimes g_+\cr
  \Delta(g_-)                   &=g_-\otimes g_++g_+\otimes g_-\cr.}$$
So $g_+$ is nonzero, grouplike in the ordinary sense, and $g_-g_+\inver$ is an odd primitive element $v_-$. So $g_+$ writes~:
$$g_+=e^{.v_+}$$
with $v_+\in V$ even. To see this one can write $g_+=1+\varepsilon$ with $\varepsilon\in S(V)\widehat\otimes\g m$ and directly check that its logarithm is primitive. We have then~:
$$g=(1+v_-)e^{.v_+}=e^{.v_-}e^{.v_+}=e^{v_-+v_+}.$$
The converse is straightforward.
\qed
\qquad Let $\widetilde{\Cal C}=S^+(V)=\bigoplus_{n\ge 1}S^n(V)$ be the cofree cocommutative graded coalgebra without co-unity cogenerated by $V$. The coproduct is given by~:
$$\widetilde\Delta(v)=\Delta(v)-1\otimes v-v\otimes 1.$$
Let $\widetilde{\Cal C}_{\sg m}=\widetilde{\Cal C}\widehat\otimes \g m$. It is easy (and left to the reader) to derive a version of the result above in that setting~:
\prop{I.2}
Any super-grouplike element in the coalgebra $\widetilde{\Cal C}_{\sg m}$ is of the form~:
$$g=e^{.v}-1=v+\cdots +{1\over n!}(v\cdots v)_{\hbox{\sevenrm n times}}+\cdots$$
with $v\in V\widehat\otimes \g m$,
and conversely any such element is super-grouplike in $\widetilde{\Cal C}_{\sg m}$.
\ndem
\qquad Let us now compute the image of a super-grouplike element by a certain coderivation in the co-unityless setting~:
\lemme{I.3}
Let $Q$ be a coderivation of coalgebra $\Cal C$ with vanishing Taylor coefficients (cf. \S\ II below) except $Q_2$, and extend it naturally by $\g m$-linearity to a coderivation of coalgebra $\widetilde C_{\sg m}$. Let $X,Y\in V\widehat\otimes \g m$ with $X$ even and $Y$ odd. Then we have~:
$$Q(e^{.(X+Y)}-1)={1\over 2}Q_2(X.X)e^{.(X+Y)}+Q_2(Y.X)e^{.X}.$$
\dem
We have the following explicit formula for a coderivation in terms of its Taylor coefficients \cite {AMM \S~III.2}~: namely, for any $n$-uple of homogeneous elements in $V$,
$$Q(x_1\cdots x_n)=\sum_{{I\coprod J=\{1,\ldots,n\}}\atop I,J\not =\emptyset}
\varepsilon_x(I,J)\bigl(Q_{|I|}(x_I)\bigr).x_J,$$
where $x_I$ stands for $x_{i_1}\cdots x_{i_k}$ when $I={\uple ik}$, and $\varepsilon_x(I,J)$ is the {\sl Quillen sign\/} associated with partition $(I,J)$, i.e. the signature of the trace on odd $x_i$'s of the shuffle associated with partition $(I,J)$. We have then~:
$$\eqalign{Q(e^{.X}-1)  &={1\over 2}Q_2(X.X)e^{.X}      \cr
        Q(Y.e^{.X})     &={1\over 2}Q_2(X.X)Y.e^{.X}+Q_2(Y.X)e^{.X}.\cr}$$
Summing up the two equalities above we get the result.
\qed
\paragraphe{II. Kontsevich's formality theorem}
Let $M$ be any real $C^\infty$ manifold, let $\g g_1$ the differential graded Lie algebra of polyvector fields on $M$ with zero differential and Schouten-Nijenhuis bracket, and let $\g g_2$ the differential graded Lie algebra of polydifferential operators on $M$ with Gerstenhaber bracket and Hochschild diffenrential. The gradings are such that a degree $n$ homogeneous element in $\g g_1$ (resp. $\g g_2$) is a $n+1$-vector field (resp. a $n+1$-differential operator).
\ssq
We can consider the shifted spaces $\g g_1[1]$ and $\g g_2[1]$ as formal graded pointed manifolds~: it means that for $i=1,2$ we have a coderivation $Q^i$ of degree $1$ on coalgebra without co-unity $S^+(\g g_i[1])$ satisfying the {\sl master equation\/}~:
$$[Q^i,Q^i]=0,$$
where $\g g_i[1]$ is meant for space $\g g_i$ with grading shifted by $1$~: a degree $n$ homogeneous element in $\g g_1[1]$ (resp. $\g g_2[1]$) is now a $n+2$-vector field (resp. a $n+2$-differential operator).
\th{II.1 (M. Kontsevich)}
There exists a $L_\infty$-quasi-isomorphism from formal graded pointed manifold $\g g_1[1]$ to formal graded pointed manifold $\g g_2[1]$~: namely, there exists a coalgebra morphism~:
$$\Cal U:S^+(\g g_1[1])\fleche 8 S^+(\g g_2[1])$$
such that~:
$$\Cal U\circ Q^1=Q^2\circ U,$$
and such that the restriction $\Cal U_1$ of $\Cal U$ to $\g g_1[1]$ is a quasi-isomorphism of cochain complexes\footnote * {\eightrm According to \cite {AMM} one should replace Schouten bracket with minus the bracket taken in the reverse order. This bracket coincides with Schouten bracket modulo a minus sign when tho odd elements are involved so it does not matter in what follows.}.
\ndem
Let us briefly recall how formality theorem is related to deformation quantization~: due to the universal property of cofree cocommutative coalgebras, coderivations $Q^i$ and $L_\infty$-quasi-isomorphism $\Cal U$ are uniquely determined by their {\sl Taylor coefficients}~:
$$\eqalign{Q_k^i:S^k(\g g_i[1])         &\fleche 8 \g g_i[2]    \cr
                \Cal U_k:S^k(\g g_1[1]) &\fleche 8 \g g_2[1],\cr}$$
$k\ge 1, i=1,2,$ obtained by composing $Q^i$ and $\Cal U$ on the right by the canonical projection~: $S^+(\g g_i)\surj 5 \g g_i$ (resp.$S^+(\g g_2)\surj 5 \g g_2$ ). Let $\g m=\hbar\R[[\hbar]]$ the projective limit of finite-dimensional nilpotent commutative algebras $\g m_r=\hbar\R[[\hbar]]/\hbar^r\R[[\hbar]]$. Let $\hbar\gamma=\hbar(\gamma_0+\hbar\gamma_1+\hbar^2\gamma_2+\cdots)$ be an infinitesimal formal Poisson bivector field, i.e. a solution in $\g g_1^{(1)}\widehat\otimes \g m$ of Maurer-Cartan equation~:
$$(\hbar d\gamma+)-{1\over 2}[\hbar\gamma,\hbar\gamma]=0,$$
which amounts exactly to the more geometrical equation~:
$$Q^1(e^{.\hbar\gamma}-1)=0,$$
where $e^{.\hbar\gamma}-1$ is grouplike (in the usual sense) in coalgebra $S^+(\g g_1[1])\widehat\otimes \g m$. Then $\Cal U(e^{.\hbar\gamma}-1)$ is grouplike in coalgebra $S^+(\g g_2[1])\widehat\otimes \g m$. So we have~:
$$ \Cal U(e^{.\hbar\gamma}-1)=e^{.\hbar\tilde \gamma}-1$$
with~:
$$\tilde \gamma=\sum_{k\ge 1}{\hbar^k\over k!}\Cal U_k(\gamma^{.k}).$$
Due to the fact that $Q_2$ vanishes at $e^{.\hbar\tilde\gamma}-1$ the element $\hbar\tilde \gamma$ verifies Maurer-Cartan equation in $\g g_2\widehat\otimes \g m$~:
$$\hbar d\tilde\gamma -{1\over 2}[\hbar\tilde \gamma,\hbar\tilde \gamma]=0.$$
We denote by $m$ the particular bidifferential  operator~:$f\otimes g\mapsto fg$, and we set $*=m+\hbar\tilde \gamma$. Maurer-Cartan equation for $\hbar\tilde \gamma$ is equivalent to~:
$$[*,*]=0,$$
i.e. $*$ is an associative product on $C^\infty(M)[[\hbar]]$. Starting from a Poisson bivector field $\gamma=P$ on $M$ we construct then explicitly a star product from $P$ and $L_\infty$-morphism $\Cal U$.
\smallskip
{\sl Remark\/}~: the expression $e^{.\hbar\gamma}-1$ is nothing but an algebraic way to express ``the point $\hbar \gamma$ in the formal graded pointed manifold''. One can be convinced by looking at the delta distribution at $\hbar \gamma$ and expressing it at $0$ by means of Taylor expansion. the expression $e^{.\hbar\gamma}-1$ is then just the difference between the delta distribution at $\hbar \gamma$ qnd the delta distribution at $0$. I would like to thank Siddhartha Sahi for having brought this nice geometrical interpretation to my attention.
\paragraphe{III. On particular super-grouplike elements}
\qquad Let $\gamma$ be a formal Poisson $2$-tensor on manifold $M$, and let $*=m+\hbar\tilde \gamma$ the star-product constructed from these data with Kontsevich's $L_\infty$-quasi-isomorphism $\Cal U$ following the formula recalled in previous paragraph. Let us consider for any $g\in C^\infty(M)[[\hbar]]$ and for any formal vector field $Y$ the super-grouplike element $e^{.\hbar(\gamma+Y+g)}-1$. we will denote by $H_g$ the hamiltonian formal vector field $[\gamma,g]$. As a straightforward application of lemma I.3 we get the following result~:
\lemme{III.1}
With the same notations as in \S\ II we have~:
$$Q^1(e^{.\hbar(\gamma+Y+g)}-1)=\hbar^2(H_g.e^{.\hbar(\gamma+Y+g)}+[Y,\gamma+g]e^{.\hbar(\gamma+g)}).$$
\ndem
\qquad 
Any morphism of graded coalgebras, in particular $L_\infty$-quasi-isomorphism $\Cal U$, preserves super-grouplike elements. Due to this fact and to proposition I.1 we have then~:
\prop{III.2}
There exists a formal differential operator $\Phi(Y)$ and a formal series $\Phi(g)\in C^\infty(M)[[\hbar]]$ such that~:
$$\Cal U(e^{.\hbar(\gamma+Y+g)}-1)=e^{.*-m+\hbar\Phi(Y)+\hbar\Phi(g)}-1,$$
with~:
$$\eqalign{\Phi(Y)      &=\Cal U_1(Y)+\hbar\Cal U_2(Y.\gamma)+{\hbar^2\over 2}\Cal U_3(Y.\gamma.\gamma)+\cdots        \cr
        \Phi(g)         &=\Cal U_1(g)+\hbar\Cal U_2(g.\gamma)+{\hbar^2\over 2}\Cal U_3(g.\gamma.\gamma)+\cdots        \cr}$$
\ndem
Correspondence $\Phi$ is precisely the tangent map at $\hbar\gamma$ of quasi-isomorphism $\Cal U$ \cite {K\S ~8.1}.
\ssq
We will now compute both terms $\pi\Cal U Q^1(e^{.\hbar(\gamma+Y+g)}-1)$ and $\pi Q^2\Cal U(e^{.\hbar(\gamma+Y+g)}-1)$, and try to get some information from the fact that they coincide, by the very definition of a $L_\infty$-morphism. Let us introduce for any pair $(Y,Z)$ of polyvector fields the second derivative term~:
$$\Psi(Y,Z)=\sum_{k\ge 0}{\hbar^k\over k!}\Cal U_{k+2}(Y.Z.\gamma^{.k}).$$
This expression is symmetric (in the graded sense) in $Y,Z\in\g g_1[1]$ and is of degree $|Y|+|Z|-2$ in $\g g_2[1]$, so it belongs to $C^\infty(M)[[\hbar]]$ when $|Y|+|Z|=0$ in $\g g_1$. The expression is skew-symmetric in $(Y,Z)$ when $Y$ and $Z$ are both vector fields, and symmetric when $Y$ is a function and $Z$ is a bivector field. We easily compute~:
$$\eqalign{
\Cal UQ^1(e^{.\hbar(\gamma+Y+g)}-1)  &=
        \hbar^2\Cal U(H_g.e^{.\hbar(\gamma+Y+g)}+[Y,\gamma+ g]e^{.\hbar(\gamma+g)})       \cr
                                &=
        \hbar^2\Cal U\bigl((H_g+\hbar H_g.Y+[Y,g]+[Y,\gamma])e^{.\hbar(\gamma+g)}\bigr).\cr}$$
From degree considerations we easily derive~:
$$\pi\Cal UQ^1(e^{.\hbar(\gamma+Y+g)}-1)=
        \hbar^2\pi\Cal U\bigl((H_g+\hbar H_g.Y+[Y,g]+[Y,\gamma]+\hbar[Y,\gamma].g)e^{.\hbar \gamma}\bigr),$$
so that we finally get~:
$$\pi\Cal UQ^1(e^{.\hbar(\gamma+Y+g)}-1)=
        \hbar^2\bigl(\Phi(H_g)+\Phi([Y,g])+\hbar \Psi(H_g,Y)
        +\Phi([Y,\gamma])+\hbar \Psi([Y,\gamma].g) \bigr).$$
On the other hand we have to compute~:
$$\eqalign{\pi Q^2\Cal U(e^{.\hbar(\gamma+Y+g)}-1)    &=\pi Q^2(e^\delta -1)  \cr
               &=[\delta, m]+{1\over 2}Q_2^2(\delta .\delta),}$$
with $\delta=*-m+\hbar(\Phi(Y)+\Phi(g))$, according to proposition III.2. We have then~:
$$\eqalign{\hbox to -12mm{}\pi Q^2\Cal U(e^{.\hbar(\gamma+Y+g)}-1)
        &=[*-m+\hbar\Phi(Y)+\hbar\Phi(g), m]\cr
        &\hbox to 12mm{}+{1\over 2}
                Q_2^2\Bigl(\bigl(*-m+\hbar\Phi(Y)+\hbar\Phi(g)\bigr).\bigl(*-m+\hbar\Phi(Y)+\hbar\Phi(g)\bigr)\Bigr)    \cr
                &\cr
        &=[*,m]+\hbar[\Phi(Y),m]+\hbar[\Phi(g),m]       \cr
        &\ \ -\dm [*,m]-\hdm [*,\Phi(Y)]+\hdm [*,\Phi(g)]       \cr
        &\ \ -\dm[m,*]+\hdm[m,\Phi(Y)]-\hdm [m,\Phi(g)] \cr
 &\ \ +\hdm [\Phi(Y),*]-\hdm [\Phi(Y),m]+{\hbar^2\over 2}[\Phi(Y),\Phi(g)] \cr
 &\ \ +\hdm [\Phi(g),*]-\hdm [\Phi(g),m]-{\hbar^2\over 2}[\Phi(Y),\Phi(g)] \cr
                &\cr
 &=-\hbar[*,\Phi(Y)]+\hbar[*,\Phi(g)]+\hbar^2[\Phi(Y),\Phi(g)].}$$
In the computation above the relation between the second Taylor coefficient $Q_2^2$ and Gerstenhaber bracket is the following~:
$$Q_2^2(x.y)=(-1)^{|x|(|y|-1)}[x,y].$$
This extra sign one must take care of comes from the identification of $S^k(\g g_2[1])$ with $\Lambda^k(\g g_2)[k]$ which goes as follows~:
$$x_1\cdots x_k\longmapsto \varepsilon.x_1\wedge\cdots\wedge x_k,$$
where $\varepsilon$ is the signature of the unshuffle storing even elements on the left and odd elements on the right \cite{AMM \S~II.4}, \cite{K \S\ 4.2}.
\ssq
We now identify the homogeneous components of degrees $1$, $0$ and $-1$ in the two expressions, so we get the following three equations~:
\th{III.3}
1) $[*,\Phi(Y)]=\hbar\Phi([\gamma,Y])$
\smallskip
2) $[*,\Phi(g)]=\hbar \Phi(H_g)$
\smallskip
3) $[\Phi(Y),\Phi(g)]=\Phi([Y,g])+\hbar\bigl(\Psi(H_g,Y)-\Psi([\gamma,Y],g)\bigr).$
\ndem
Equation 1) implies the following : tangent map $\Phi$ sends derivations of $C^\infty(M)[[\hbar]]$ with commutative product leaving $\gamma$ invariant to derivations of $C^\infty(M)[[\hbar]]$ with deformed product. From equation 2) we see that hamiltonian formal vector fields are sent to inner derivations of the deformed algebra~: these two facts proceed more directly from the fact that the tangent map is a morphism of cochain complexes \cite {K \S\ 8.1}. Equation 3) rewrites as follows~:
$$R(Y,g)={1\over \hbar}\bigl(\Phi([Y,g])-[\Phi(Y),\Phi(g)]\bigr),$$
where skew-symmetric bilinear term~:
$$R(Y,g)=\Psi([\gamma,g],Y)-\Psi([\gamma,Y],g)$$
can be seen as a kind of curvature.  As a particular case we can take for $Y$ the hamiltonian vector field $H_f$ for any $f\in C^\infty(M)[[\hbar]]$. From equations 2) and 3) we immediately get~:
\th{III.4}
For any $f,g$ in $C^\infty(M)$ the following formula holds true~:
$$\Psi(H_f,H_g)={1\over \hbar}\Bigl(\Phi(\{f,g\})-{\Phi(f)*\Phi(g)-\Phi(g)*\Phi(f)\over\hbar}\Bigr).$$
\ndem
{\sl Remark 1\/}~: Introducing the new star-product~:
$$f\#g=\Phi\inver\bigl(\Phi(f)*\Phi(g) \bigr)$$
the formula of theorem III.4 can be rewritten as follows~:
$$\Phi\inver\circ \Psi(H_f,H_g)={1\over \hbar}\Bigl(\{f,g\}-{f\#g-g\#f\over \hbar} \Bigr).$$
{\sl Remark 2\/}~: It is immediate to see from symmetry properties of star products that the right-hand side is $O(\hbar)$. In the flat case $M=\R^d$ one can use quasi-isomorphism of \cite {K \S ~6}, and recover this fact from the left-hand side by computing the constant term~: it involves only the following graph~:
\diagramme
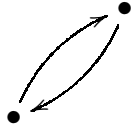
{
\xymatrix{
& \bullet \ar @/^/[dl]\\
\bullet \ar @/^/[ur]}
}
the weight of which is zero \cite {K \S ~7.3.1}.
\paragraphe{IV. Gauge transformations}
Recall (from \cite {K \S\ 3.2}) that the {\sl gauge group\/} associated with any differential graded Lie algebra $\g g$ is by definition the pro-nilpotent group $G$ associated with pro-nilpotent Lie algebra $\g g^{(0)}\widehat\otimes \g m$. The gauge group acts on $\g g^{(1)}\widehat\otimes \g m$ by affine transformations, and the action of $G$ is defined by exponentiation of the infinitesimal action of $\g g^{(0)}\widehat\otimes \g m$~:
$$\alpha\otimes\gamma\in\g g^{(0)}\otimes\g g^{(1)}\mapsto
        \alpha.\gamma:=d\alpha+[\alpha,\gamma].$$
It is easy to check that gauge group $G$ acts on the set of solutions of Maurer-Cartan equation in $\g g^{(1)}\widehat\otimes \g m$.    
\ssq
We try in this last paragraph to give a geometrical meaning to equation 1) of theorem III.3. We let $\hbar\gamma$ run inside the set $(MC)_1$ of infinitesimal formal Poisson $2$-tensors on manifold $M$, i.e. the subset of solutions of Maurer-Cartan equation in $\g g_1^{(1)}\widehat\otimes \g m$. We denote by $*_\gamma$ the star-product constructed from $\gamma$ along the lines of \S\ II. We will use the notation $\Cal U'_{\hbar\gamma}$ instead of $\Phi$ to emphasize the dependance on $\gamma$. Equation 1) of theorem III.3 is then rewritten as follows~:
$$[\Cal U'_{\hbar\gamma}(Y),*_\gamma]=\hbar\Cal U'_{\hbar\gamma}([Y,\gamma])$$
for any vector field $Y$ on manifold $M$. Let $\Cal V_Y^1$ be the vector field (of degree $0$) on formal pointed manifold $\g g_1^{(1)}\widehat\otimes \g m$ equal to $[Y,\hbar\gamma]$ at point $\hbar \gamma$. It is the coderivation of coalgebra $S^+(\g g_1^{(1)})\widehat\otimes \g m$ such that~:
$$\Cal V_Y^1(e^{.\hbar\gamma}-1)=[Y,\hbar\gamma]e^{.\hbar\gamma}.$$
It restricts to the submanifold $(MC)_1$.
\ssq
The $L_\infty$-quasi-isomorphism $\Cal U$ is injective, by injectivity of first Taylor coefficient $\Cal U_1$. Let $(MC)_2$ be the set of solutions of Maurer-Cartan equation in $\g g_2^{(1)}\widehat\otimes \g m$. Let $\Cal V_Y^2$ be the vector field (of degree $0$) on formal pointed manifold $\Cal U(\g g_1^{(1)}\widehat\otimes \g m)$ equal to $[\Cal U'_{\hbar\gamma}(Y), *_\gamma]$ at point $\hbar\tilde\gamma$. It is the coderivation of the image coalgebra $\Cal U\bigl(S^+(\g g_1^{(1)})\widehat\otimes \g m \bigr)$ such that~:
$$\Cal V_Y^2(e^{.\hbar\tilde\gamma}-1)=[\Cal U'_{\hbar\gamma}(Y),*_\gamma]e^{.\hbar\tilde\gamma}.$$
Clearly vector field $\Cal V_Y^2$ restricts to $\Cal U\bigl((MC)_1\bigr)\subset (MC)_2$. Adding multiplication $m$ we identify $(MC)_2$ with the set of all star products, and $\Cal U\bigl((MC)_1\bigr)$ with the set $(FSP)$ of formality star products. 
\prop{IV.1}
$$\Cal U\circ \Cal V_Y^1=\Cal V_Y^2\circ \Cal U.$$
\dem
We have~:
$$\Cal U\circ \Cal V_Y^1 (e^{.\hbar\gamma}-1)=\Cal U'_\gamma([Y,\hbar\gamma])
                                        e^{.\hbar\tilde\gamma}$$
and~:
$$\Cal V_Y^2\circ \Cal U(e^{.\hbar\gamma}-1)=[\Cal U'_{\hbar\gamma}(Y),*_\gamma]
                                       e^{.\hbar\tilde\gamma}.$$
The result follows then immediately from equation 1) of theorem III.3.
\qed
It is clear that we have $\Cal V_Y^2(\tau)=O(\hbar)$ for any $\tau\in \Cal U \bigl((MC)_1\bigr)$. We have then~:
\cor{IV.2}
Let $\gamma_Y=e^{\smop{ad}\hbar Y}\gamma$ be the transformation of formal Poisson $2$-tensor $\gamma$ under the formal diffeomorphism $e^{\hbar Y}$. It clearly belongs to $(MC)_1$ and the star-product $*_{\gamma_Y}$ constructed from $\gamma_Y$ by means of $L_\infty$-morphism $\Cal U$ is the transformation of star-product $*_\gamma$ under the formal diffeomorphism $e^{\Cal V_Y^2}$ of formal pointed manifold $(FSP)$.
\ndem
\qquad
Of course differential operator $\Cal U'_\gamma(Y)$ is not a vector field on $M$ in general, so diffeomorphism $e^{\Cal V_Y^2}$ of formal pointed manifold $(FSP)$ does not come from a formal diffeomorphism of $M$. Correspondences $Y\mapsto \Cal V_Y^1$ ans $Y\mapsto \Cal V_Y^2$ are injective and respect brackets, i.e.:
$$\Cal V_{[Y,Z]}^1=[\Cal V_Y^1,\Cal V_Z^1],$$
and, as a consequence of proposition IV.1 :
$$\Cal V_{[Y,Z]}^2=[\Cal V_Y^2,\Cal V_Z^2].$$
By exponentiation we get then the following result~:
\th{IV.3}
Let $G_1$ and $G_2$ denote the gauge groups of $\g g_1$ and $\g g_2$ respectively, and let $\Cal G$ be the group of formal diffeomorphisms of $(FSP)$. The correspondence~:
$$\eqalign{\iota:G_1    &\longrightarrow \Cal G       \cr
        e^Y             &\longmapsto e^{\Cal V_Y^2}     \cr}$$ 
is an embedding of $G_1$ into the group $\Cal G$ of formal diffeomorphisms of $(FSP)$ such that~:
$$*_{g.\gamma}=\iota(g).*_\gamma.$$
\ndem
{\sl Remark\/}~: On one hand we have natural embeddings~:
$$G_1\subset G_2,\hbox to 12mm{} G_1\subset \Cal G,$$
but we must stress on the other hand that vector fields $\Cal V_Y^1$ are linear whereas vector fields $\Cal V_Y^2$ are not, and not even affine (this is due to the non-linearity of $\Cal U$). The embedding of $G_1$ constructed above is then nontrivial, in the sense that the image of $G_1$ is not even contained in the second gauge group $G_2$.

\vskip 12mm
\paragraphe{References}
\bib{AMM} D. Arnal, D. Manchon, M. Masmoudi~: {\sl Choix des signes dans la formalit\'e de Kontsevich}, eprint math.QA/0003003.
\bib{BFFLS} F. Bayen, M. Flato, C. Fr\o nsdal, A. Lichnerowicz, D. Sternheimer, {\sl Deformation theory and quantization I. Deformations of symplectic structures}, Ann. Phys. 111 No 1, 61-110 (1978).
\bib{K} M. Kontsevich, {\sl Deformation quantization of Poisson manifolds I}, eprint math.QA.9709040.

\bye